\documentclass[10pt,a4paper]{article}

\usepackage{isorot}

\usepackage{graphicx,subfigure}

\usepackage[square,numbers]{natbib}

\usepackage{booktabs}
\usepackage{multirow}

\usepackage{pstricks,pst-plot}

\usepackage{amsmath}
\usepackage{amsthm}
\usepackage{amsfonts}
\usepackage{amssymb}
\usepackage{ipa}
\usepackage{latexsym}
\usepackage{stackrel}
\usepackage[T1]{fontenc}
\usepackage{IEEEtrantools}
\usepackage{enumerate}
\usepackage{fourier}
\usepackage{skak}

\usepackage{lastpage}

\let\mathnumsetfont\mathbb
\newcommand\Nset{\mathnumsetfont N} 
\newcommand\Rset{\mathnumsetfont R} 

\newcommand\M{{\cal M}}

\newcommand\thelr{the L'H\^{o}pital's rule}

\newtheorem{thm}{Theorem}
\newtheorem*{thmck*}{Theorem CK}
\newtheorem*{thmckpm*}{Theorem CK$_{\pm\infty}$}
\newtheorem{teo*}{Theorem}
\newtheorem{prop}{Proposition}
\newtheorem{cor}{Corollary}
\newtheorem{rmq}{Remark}

\newtheorem{exm}{\emph{Example}}

\def\limsup{\mathop{\overline{\mathrm{lim}}}}
\def\liminf{\mathop{\underline{\mathrm{lim}}}}
\def\limx{\lim_{x\rightarrow\infty}}
\def\limy{\lim_{y\rightarrow\infty}}

\def\limt{\lim_{t\rightarrow\infty}}
\def\lims{\lim_{s\rightarrow\infty}}

\def\limsupx{\limsup_{x\rightarrow\infty}}
\def\liminfx{\liminf_{x\rightarrow\infty}}

\def\limsupt{\limsup_{t\rightarrow\infty}}
\def\liminft{\liminf_{t\rightarrow\infty}}

\def\limn{\lim_{n\rightarrow\infty}}
\def\limsupn{\limsup_{n\rightarrow\infty}}

\def\limk{\lim_{k\rightarrow\infty}}

\def\barF{\overline{F}}

\title{Revisiting extensions of regularly varying functions}

\author{Meitner Cadena\thanks{UPMC Paris 6 \& CREAR, ESSEC Business School;\, E-mail: meitner.cadena@etu.upmc.fr, b00454799@essec.edu, or meitner.cadena@gmail.com}} 

\date{}

\begin{document}

\maketitle

\begin{abstract}

Relationships among the classes ${\cal M}$, ${\cal M}_\infty$, and ${\cal M}_{-\infty}$ and the class of \emph{O}-regularly varying functions are shown.
These results are based on two characterizations of ${\cal M}$, ${\cal M}_\infty$, and ${\cal M}_{-\infty}$ provided by Cadena and Kratz in \cite{NN2014} and a new one given in this note.

\vspace{2ex}

{\it Keywords: regularly varying function; slowly varying function; O-RV; large deviations; extreme value theory}

\vspace{2ex}

{\it AMS classification}:  26A12; 60F10
\end{abstract}

\section{Introduction}

A positive and measurable function $U$ defined on $\Rset^+$ is a \emph{regularly varying} (RV) function if
\begin{equation}\label{eq:20150102:01}
\limx\frac{U(tx)}{U(x)}<\infty\textrm{\quad($t>0$).}
\end{equation}
If this limit equals 1, $U$ is a \emph{slowly varying} (SW) function.
Classes \emph{RV} and \emph{SV} of regularly and slowly varying functions were introduced by Karamata \cite{Karamata1930} in 1930.
Since then theory of these functions has been developed in many directions.
Systematic treatment of this theory can be found in e.g. \cite{BinghamGoldieTeugels} and \cite{Seneta}.

Extensions of RV functions have been obtained by letting \eqref{eq:20150102:01} to vary.
An early extension of this type was given by Avakumovi\'{c} in 1936 \cite{Avakumovic1936}.
He introduced the class \emph{O-RV} of \emph{O-regularly varying} (O-RV) functions $U$ which satisfy the following condition instead of \eqref{eq:20150102:01}:
\begin{equation}\label{eq:20150102:02}
0<U_*(t):=\liminfx\frac{U(tx)}{U(x)}\leq\limsupx\frac{U(tx)}{U(x)}=:U^*(t)<\infty\textrm{\quad($t\geq1$).}
\end{equation}
Recently Cadena and Kratz \cite{NN2014} gave an extension of RV functions by also letting \eqref{eq:20150102:01} to vary,
but they designed it in a different way to the previous one.
They introduced the class $\M$ which consists in functions $U$ satisfying the following condition instead of \eqref{eq:20150102:01}:
\begin{equation}\label{eq:20150102:03}
\exists\rho\in\Rset\textrm{, }\forall\epsilon>0\textrm{, }\limx\frac{U(x)}{x^{\rho+\epsilon}}=0\quad\textrm{and}\quad\limx\frac{U(x)}{x^{\rho-\epsilon}}=\infty\textrm{.}
\end{equation}
We have clearly $\textrm{\emph{RV}}\subsetneq\textrm{\emph{O-RV}}$ and, for instance using Theorem \ref{teo:20141211:001} (see Corollary \ref{cor:20141225:001}),
$\textrm{RV}\subsetneq\M$.
There arises the natural question of how \emph{O-RV} and $\M$ are related between them.
We undertake this study helping us of characterizations of these classes: recalling well-known characterizations of \emph{O-RV} and giving proofs of three characterizations of $\M$,
two of them provided in \cite{NN2014} and a new one given in this note.

Cadena and Kratz also introduced the following natural extensions of $\M$.
\begin{eqnarray}
\M_{\infty} & := & \left\{U:\Rset^+\to\Rset^+:\textrm{$U$ is measurable and satisfies }\forall\rho\in\Rset\textrm{, }\limx\frac{U(x)}{x^{\rho}}=0\right\} \label{eq:20150102:04} \\
\M_{-\infty} & := & \left\{U:\Rset^+\to\Rset^+:\textrm{$U$ is measurable and satisfies }\forall\rho\in\Rset\textrm{, }\limx\frac{U(x)}{x^{\rho}}=\infty\right\} \label{eq:20150102:05}\textrm{.}
\end{eqnarray}
The new characterization given for $\M$ is extended to $\M_\infty$ and $\M_{-\infty}$.
Relationships among $\M_\infty$ and $\M_{-\infty}$ and \emph{O-RV} are also investigated in this note.

This note is organized as follows. 
The main results are presented in the next section, introducing previously notations and definitions.
First, the new characterizations of $\M$, $\M_\infty$, and $\M_{-\infty}$ based on limits are given.
Next, analyses of uniform convergence in these characterizations are presented and, finally, relationships among \emph{O-RV} and $\M$, $\M_\infty$ and $\M_{-\infty}$ are shown.
All proofs are collected in Section \ref{SectionProofs}.
Conclusion is presented in the last section.

\section{Main Results}

For a positive function $U$ with support $\Rset^+$ its \emph{lower} and \emph{upper orders} are defined by (see e.g. \cite{BinghamGoldieTeugels})
$$
\begin{array}{rllrll}
\mu(U) & := & \displaystyle \liminfx\frac{\log\left(U(x)\right)}{\log(x)}\textrm{,} & 
\nu(U) & := & \displaystyle \limsupx\frac{\log\left(U(x)\right)}{\log(x)}\textrm{.}
\end{array}
$$
Throughout this note $\log(x)$ represents the natural logarithm of $x$.

We notice that the classes $\M$, $\M_\infty$, and $\M_{-\infty}$ defined in \eqref{eq:20150102:03}, \eqref{eq:20150102:04}, and \eqref{eq:20150102:05} are a bit weaker than the corresponding classes given in \cite{NN2014}, and that each of them is disjoint from each other.
Moreover, using straightforward computations, one can prove that $\rho$ defined in \eqref{eq:20150102:03} is unique, hence it will be denoted by $\rho_U$, and one can show that
$\epsilon>0$ in \eqref{eq:20150102:03} can be taken sufficiently small.
Additionally, one can prove that $\M$ is strictly larger than RV, for instance using Theorem \ref{teo:20141211:001} (see Corollary \ref{cor:20141225:001}),
and that $\M_\infty$ is related to the domain of attraction of Gumbel (see \cite{NN2014}).

The new characterizations of $\M$, $\M_\infty$, and $\M_{-\infty}$ follow.

\begin{thm}\label{teo:20141211:001}
Let $U:\Rset^+\to\Rset^+$ be a measurable function. Then
\begin{itemize}
\item[(i)]
$U\in\M$ with $\rho_U=-\tau$ iff
\begin{equation}\label{eq:20141213:001}
\left\{
\begin{array}{cl}
 & \displaystyle \forall r<\tau\textrm{, }\exists x_a>1\textrm{, }\forall x\geq x_a\textrm{, }\limt t^r\frac{U(tx)}{U(x)}=0 \\
\textrm{} & \\
 & \displaystyle \forall r>\tau\textrm{, }\exists x_b>1\textrm{, }\forall x\geq x_b\textrm{, }\limt t^r\frac{U(tx)}{U(x)}=\infty\textrm{.}
\end{array}
\right.
\end{equation}

\item[(ii)]
$U\in\M_\infty$ iff
\begin{equation}\label{eq:20141213:002}
\forall r\in\Rset\textrm{, }\exists x_0>1\textrm{, }\forall x\geq x_0\textrm{, }\limt t^r\frac{U(tx)}{U(x)}=0\textrm{.}
\end{equation}

\item[(iii)]
$U\in\M_{-\infty}$ iff \ 
\begin{equation}\label{eq:20141213:003}
\forall r\in\Rset\textrm{, }\exists x_0>1\textrm{, }\forall x\geq x_0\textrm{, }\limt t^r\frac{U(tx)}{U(x)}=\infty\textrm{.}
\end{equation}

\end{itemize}
\end{thm}

\begin{exm}\label{ex:20141207:001}~

\begin{enumerate}
\item
Consider a measurable and positive function $U$ with support $\Rset^+$ such that, for $x\geq x_0$ with some $x_0>1$, $U(x)=x\big/\log(x)$.

Noting that, for $t,x>1$,
$$
\limt t^r\frac{U(tx)}{U(x)}=
\limt t^{r-1}\frac{\log(x)}{\log(tx)}=
\left\{
\begin{array}{ll}
0 & \textrm{if $r>1$} \\
\infty & \textrm{if $r<1$,}
\end{array}
\right.
$$
provides, taking $\tau=-1$ and applying Theorem \ref{teo:20141211:001}, (i), $U\in\M$ with $\rho_U=1$.

\item
Let $U$ be a function defined by $U(x):=x^{\sin(x)}$, $x>0$.

Writing
$$
t^r\frac{U(tx)}{U(x)}=
t^{r+\sin(tx)}x^{\sin(tx)-\sin(x)}
$$
gives, for $r\in\Rset$,
$$
\limsupt t^r\frac{U(tx)}{U(x)}=\infty
\quad\textrm{and}\quad
\liminft t^r\frac{U(tx)}{U(x)}=0\textrm{.}
$$
Hence the necessary condition of Theorem \ref{teo:20141211:001}, (i), is not satisfied and consequently gives $U\not\in\M$.

\end{enumerate}
\end{exm}

It follows a consequence of Theorem \ref{teo:20141211:001}.
This result was proved by Cadena and Kratz in \cite{NN2014} combining a result provided in \cite{deHaan} and another characterization of $\M$ (see Theorem CK later).

\begin{cor}\label{cor:20141225:001}
RV \ $\subsetneq$ \ $\M$.
\end{cor}

Note that, from Corollary \ref{cor:20141225:001}, $\textrm{\emph{RV}}\subseteq\M\bigcap\textrm{\emph{O-RV}}$.

There are not common elements between \emph{O-RV} and $\M$ under their definitions given in \eqref{eq:20150102:02} and \eqref{eq:20150102:03} respectively,
but observing the characterization of $\M$ given in Theorem \ref{teo:20141211:001} one identifies the quotient $U(tx)\big/U(x)$, which appears in \eqref{eq:20150102:02}.
The next example exploits this link to show a first relationship between \emph{O-RV} and $\M$.

\begin{exm}\label{exm:20150103:001}
$\M$ \ $\not\subseteq$ \ O-RV.

Let $U$ be a function defined by $U(x):=\exp\big\{(\log x)^{\alpha}\,\cos\big((\log x)^{\beta}\big)\big\}$, $x>0$, where $0<\alpha,\beta<1$ such that $\alpha+\beta>1$.

Prof. Philippe Soulier gave recommendations to correct an error in an early version of this example.

On the one hand, noting that, for $x,t>e$, using the changes of variable $y=\log(x)$ and $s=\log(t)$ and observing that $s\to\infty$ as $t\to\infty$,
\begin{eqnarray*}
\lefteqn{\limt t^r\frac{U(tx)}{U(x)}=
\lims \exp\left\{rs+(s+y)^{\alpha}\,\cos\big((s+y)^{\beta}\big)-y^{\alpha}\,\cos\big(y^{\beta}\big)\right\}} \\
 & & =
\lims \exp\left\{s\left(r+\frac{1}{s^{1-\alpha}}\left(1+\frac{y}{s}\right)^{\alpha}\,\cos\big((s+y)^{1/3}\big)-\frac{y^{\alpha}}{s}\,\cos\big(y^{\beta}\big)\right)\right\}=
\left\{
\begin{array}{ll}
0 & \textrm{if $r>0$} \\
\infty & \textrm{if $r<0$,}
\end{array}
\right.
\end{eqnarray*}
provides, taking $\tau=0$ and applying Theorem \ref{teo:20141211:001}, $U\in\M$ with $\rho_U=0$.

On the other hand, writing, for $x>e$ and $t>0$, using the previous changes of variables, with $x$ such that $\big(\log tx\big)^{\beta}=\pi\big/2+2k\pi$, for a given $t$,

\begin{eqnarray*}
\lefteqn{\frac{U(tx)}{U(x)} 
=\exp\left\{-\big(\log x\big)^{\alpha}\cos\big((\log x)^{\beta}\big)\right\}} \\
 & & =\exp\left\{-y^{\alpha}\cos\big(((\pi\big/2+2k\pi)^{1/\beta}-s)^\beta\big)\right\} \\
 & & =\exp\left\{\big((\pi\big/2+2k\pi)^{1/\beta}-s\big)^{\alpha}\sin\big(((\pi\big/2+2k\pi)^{1/\beta}-s)^\beta-(\pi\big/2+2k\pi)\big)\right\}\textrm{.}
\end{eqnarray*}

Since $\big((\pi\big/2+2k\pi)^{1/\beta}-s\big)^\beta-(\pi\big/2+2k\pi)\to0$ as $k\to\infty$, 
we have


\begin{eqnarray*}
\lefteqn{\limk\left[\big((\pi\big/2+2k\pi)^{1/\beta}-s\big)^{\alpha}\sin\big(((\pi\big/2+2k\pi)^{1/\beta}-s)^\beta-(\pi\big/2+2k\pi)\big)\right]} \\
 & = & 
\limk\frac{((\pi\big/2+2k\pi)^{1/\beta}-s)^\beta-(\pi\big/2+2k\pi)}{\big((\pi\big/2+2k\pi)^{1/\beta}-s\big)^{-\alpha}}\textrm{,}
\end{eqnarray*}

which is an indetermination of type $0\big/0$.
Then, applying L'Hopital's rule we have

\begin{eqnarray*}
\lefteqn{\limk\frac{((\pi\big/2+2k\pi)^{1/\beta}-s)^\beta-(\pi\big/2+2k\pi)}{\big((\pi\big/2+2k\pi)^{1/\beta}-s\big)^{-\alpha}}} \\
 & = & \limk(2\pi)^{\alpha/\beta}\frac{((\pi\big/2+2k\pi)^{1/\beta}-s)^\beta-(\pi\big/2+2k\pi)}{k^{-\alpha/\beta}} \\
 & = & \limk-\frac{\beta}{\alpha}(2\pi)^{\alpha/\beta+1}\frac{((\pi\big/2+2k\pi)^{1/\beta}-s)^{\beta-1}(\pi\big/2+2k\pi)^{1/\beta-1}-1}{k^{-\alpha/\beta-1}}\textrm{,}
\end{eqnarray*}

which is an indetermination of type $0\big/0$.
Then, applying again L'Hopital's rule we have

\begin{eqnarray*}
\lefteqn{\limk\frac{((\pi\big/2+2k\pi)^{1/\beta}-s)^\beta-(\pi\big/2+2k\pi)}{\big((\pi\big/2+2k\pi)^{1/\beta}-s\big)^{-\alpha}}} \\
& = & \limk s\frac{\beta(1-\beta)}{\alpha(\alpha+\beta)}(2\pi)^{\alpha/\beta+2}\frac{((\pi\big/2+2k\pi)^{1/\beta}-s)^{\beta-2}(\pi\big/2+2k\pi)^{1/\beta-2}}{k^{-\alpha/\beta-2}} \\
& = & \limk s\frac{\beta(1-\beta)}{\alpha(\alpha+\beta)}(2\pi)^{(\alpha+\beta-1)/\beta}k^{(\alpha+\beta-1)/\beta} \\
& = & \left\{
\begin{array}{ll}
\infty & \textrm{if s>0} \\
-\infty & \textrm{if s<0.}
\end{array}
\right.\textrm{}
\end{eqnarray*}

Then, we get, for $t>1$,
$$
U^*(t)=\limsupx\frac{U(tx)}{U(x)}=\infty\textrm{,}
$$
and, for $t<1$,
$$
U_*(t)=\liminfx\frac{U(tx)}{U(x)}=0\textrm{,}
$$
which contradict \eqref{eq:20150102:02}, so $U\not\in\textrm{O-RV}$.
In particular, $U\not\in\textrm{SV}$.
\end{exm}

Next, the uniform convergences in $x$ of limits given in \eqref{eq:20141213:001}, \eqref{eq:20141213:002}, and \eqref{eq:20141213:003} are analyzed.
To this aim, we will use the next two results.

\begin{prop}\label{prop:20141211:001}
Let $U:\Rset^+\to\Rset^+$ be a measurable function. Then
\begin{itemize}
\item[(i)]
If $U\in\M$ with $\rho_U=-\tau$, then there exists $x_0>1$ such that, for $x_0\leq c< d<\infty$, there exist $0<M_c<M_d$ satisfying, for $x\in[c;d]$, $M_c\leq U(x)\leq M_d$.

\item[(ii)]
If $U\in\M_{\infty}$, then there exists $x_0>1$ such that, for $c\geq x_0$, there exist $M_c>0$ satisfying, for $x\in[c;\infty)$, $U(x)\leq M_c$.

\item[(iii)]
If $U\in\M_{-\infty}$, then there exists $x_0>1$ such that, for $d\geq x_0$, there exist $M_d>0$ satisfying, for $x\in[d;\infty)$, $U(x)\geq M_d$.

\end{itemize}
\end{prop}

\begin{prop}[Given in \cite{Arandelovic2004}]\label{propArandelovic2004}
Let $\mu$ be the Lebesgue measure on $\Rset$, $A$ a measurable set of positive measure, and $\big\{x_n\big\}_{n\in\Nset}$ a bounded sequence of real numbers.
Then, $\mu(A)\leq\mu\big(\limsupn(x_n+A)\big)$.
\end{prop}

Now the results on uniform convergences are presented. Their proofs are inspired by \cite{ArandelovicPetkovic2007}.

\begin{thm}[Uniform Convergence Theorem (UCT)]\label{teo:20141211:002}~
Let $U:\Rset^+\to\Rset^+$ be a measurable function. Then
\begin{itemize}
\item[(i)]
If $U\in\M$ with $\rho_U=-\tau$ and $r<\tau$, then, for any $x_a\leq c< d<\infty$ for some $x_a>1$,
$$
\limt t^r\sup_{x\in[c;d]}\frac{U(tx)}{U(x)}=0\textrm{.}
$$

\item[(ii)]
If $U\in\M$ with $\rho_U=-\tau$ and $r>\tau$, then, for any $x_b\leq c< d<\infty$ for some $x_b>1$,
$$
\limt t^r\inf_{x\in[c;d]}\frac{U(tx)}{U(x)}=\infty\textrm{.}
$$

\item[(iii)]
If $U\in\M_\infty$ satisfying, for $s>1$, $U(x)\geq M_s$ for $x\in[1;s]$ and some $M_s>0$, then, for $r\in\Rset$ and any constants $x_0\leq c< d<\infty$ for some $x_0>1$,
$$
\limt t^r\sup_{x\in[c;d]}\frac{U(tx)}{U(x)}=0\textrm{.}
$$

\item[(iv)]
If $U\in\M_{-\infty}$ satisfying, for $s>1$, $U(x)\leq M_s$ for $x\in[1;s]$ and some $M_s>0$, then, for $r\in\Rset$ and any constants $x_0\leq c< d<\infty$ for some $x_0$,
$$
\limt t^r\inf_{x\in[c;d]}\frac{U(tx)}{U(x)}=\infty\textrm{.}
$$

\end{itemize}
\end{thm}

Note that UCT cannot be extended to infinite intervals.
For instance, from the function $U$ given in Example \ref{exm:20150103:001} we have that computing the supremum of the quotient $U(tx)\big/U(x)$ in $x$ on $[x_0;\infty)$, for any $x_0>1$, gives always $\infty$, and hence one cannot deduce that $\rho_U=0$.

The next results on \emph{O-RV}, $\M$, $\M_\infty$, and $\M_{-\infty}$ will be used to give more relationships between these classes.
On \emph{O-RV} we need:

\begin{prop}[see e.g. \cite{Karamata1936}, \cite{Seneta}, \cite{AljancicArandelovic1977}, \cite{GelukdeHaan}, and \cite{BinghamGoldieTeugels}]\label{prop:20141220:000}
Let $U:\Rset^+\to\Rset^+$ be a measurable function.
Then the following statements are equivalent:
\begin{itemize}
\item[(i)]
$U\in\textrm{O-RV}$.

\item[(ii)]
There exist $\alpha,\beta\in\Rset$ and $x_0>1,c>0$ such that, for all $t\geq1$ and $x\geq x_0$,
$$
c^{-1}t^\beta\leq\frac{U(tx)}{U(x)}\leq ct^\alpha\textrm{.}
$$

\item[(iii)]
There exist functions $\eta(x)$ and $\phi(x)$ bounded on $[x_0;\infty)$, for some $x_0\geq1$, such that
$$
U(x)=\exp\left\{\eta(x)+\int_1^x\phi(y)\frac{dy}{y}\right\}\textrm{,\quad$x\geq1$.}
$$
\end{itemize}
\end{prop}

On $\M$ we need the next two characterizations of $\M$ given by Cadena and Kratz in \cite{NN2014}.
For the sake of completeness of this note, we give them as Theorem CK and indicate their proofs.
Part of these proofs are copied from \cite{NN2014}.

\begin{thmck*}\label{teock}
Let $U:\Rset^+\to\Rset^+$ be a measurable function. Then the following statements are equivalent:
\begin{itemize}
\item[(i)]
$U\in\M$ with $\rho_U=\tau$.

\item[(ii)]
$
\displaystyle
\limx\frac{\log\left(U(x)\right)}{\log(x)}=\tau\textrm{.}
$

\item[(iii)]
There exist $b>1$ and measurable functions $\alpha$, $\beta$, and $\delta$ satisfying, as $x\to\infty$,
$$
\alpha(x)\big/\log(x)\to0\textrm{,}\quad\beta(x)\to\tau\textrm{,}\quad\delta(x)\to1\textrm{,}
$$
such that
$$
U(x)=\exp\left\{\alpha(x)+\delta(x)\int_b^x\beta(s)\frac{ds}{s}\right\}\textrm{,\quad$x\geq x_1$ for some $x_1\geq b$.}
$$

\end{itemize}
\end{thmck*}

\begin{rmq}
If $\barF$ is the tail of a distribution $F$ associated to a random variable (rv) $X$, some authors (see e.g. \cite{Nakagawa2007} and \cite{Nakagawa2008}) say that $X$ is heavy-tailed if the limit
$$
\eta:=\limx\frac{\log\left(\barF(x)\right)}{\log(x)}\textrm{}
$$
exists and takes a negative value.

We notice that this characterization does not cover rvs with heavy tails satisfying $\eta=0$ or with heavy tails for which such limit does not exist.
Indeed, on the one side, from Theorem CK one has that $\eta=0$ implies that $\barF\in\M$ with $\rho_{\barF}=0$, being a particular case of these functions the SV functions,
which are considered heavy-tailed.
On the other side, Cadena and Kratz presented in \cite{NN2014} families of tails $\barF$ for which the limit $\displaystyle \limx\frac{\log\left(\barF(x)\right)}{\log(x)}$ does not exist, for instance the next tail defined by (see \cite{NN2014})
\begin{quote}
Let $\alpha>0$, $\beta<-1$, $x_a>1$, and define the series $x_n=x_a^{(1+\alpha)^n}$, $n\geq1$, which satisfies $x_n\to\infty$ as $n\to\infty$.
It is not hard to prove that the tail $\barF$ associated to a rv $X$ and defined by
$$
\barF(x):=\left\{
\begin{array}{ll}
1 & x\in[0;x_1) \\
x_n^{\alpha(1+\beta)} & x\in[x_n;x_{n+1})\textrm{, }\forall n\geq1
\end{array}
\right.
$$
satisfies
$$
\liminfx\frac{\log\left(\barF(x)\right)}{\log(x)}=-\frac{\alpha(1+\beta)}{1+\alpha}<
-\alpha(1+\beta)=\limsupx\frac{\log\left(\barF(x)\right)}{\log(x)}\textrm{.}
$$
Note that if $-\alpha(1+\beta)\big/(1+\alpha)<1$, then the expected value of $X$ is $\infty$, which means that $X$ can be considered as a heavy-tailed rv.
\end{quote}
\end{rmq}

We notice from the representations of $U$ via \emph{O-RV} and $\M$ given in Proposition \ref{prop:20141220:000}, (iii), and Theorem CK, (iii), respectively,
that a key difference between those representations is the presence of a bounded function under the integral symbol.
Motivated by this observation, we built the next function belonging to \emph{O-RV} but not to $\M$.
This aim is reached by building a bounded function $\phi$ such that the limit
$
\displaystyle
\limx\frac{\int_1^x\phi(s)\frac{ds}{s}}{\log(x)}
$
does not exist.
Note that if this limit exists, then, applying Theorem CK, (iii), $U\in\M$.

\begin{exm}\label{prop:20150102:02}
O-RV \ $\not\subseteq$ \ $\M$.

Let $U:\Rset^+\to\Rset^+$ be a measurable function satisfying, for $x\geq1$, $\displaystyle U(x)=\exp\left\{\int_1^x\phi(s)\frac{ds}{s}\right\}$, where the function $\phi$ has support $[1;\infty)$ and is defined by
$$
\phi(x)=\left\{
\begin{array}{ll}
0 & \textrm{if $x\in[1;e)$ or $x\in I_n$ with $n$ odd} \\
1 & \textrm{if $x\in I_n$ with $n$ even,}
\end{array}
\right.
$$
where $I_n=[e^{e^n};e^{e^{n+1}})$, $n\in\Nset$.

On the one hand, applying Proposition \ref{prop:20141220:000}, one has $U\in\textrm{O-RV}$.

On the other hand, writing, for $x>1$, using the change of variable $y=\log(s)\big/\log(x)$,
$$
\frac{\log(U(x))}{\log(x)}=
\frac{\int_1^x\phi(s)\frac{ds}{s}}{\log(x)}=
\int_0^1\phi\left(e^{y\,\log(x)}\right)dy
$$
gives, taking $x_n=e^{e^{n}}$, $n=2,3,\ldots$,
$$
\frac{\log(U(x_n))}{\log(x_n)}=
\sum_{k=1}^{n-1}\int_{e^k/e^n}^{e^{k+1}/e^n}\phi\left(e^{y\,e^{n}}\right)dy
=\left\{
\begin{array}{ll}
\displaystyle \sum_{k=0}^{n-1}(-1)^ke^{-k} & \textrm{if $n$ is odd} \\
 & \\
\displaystyle \sum_{k=1}^{n-1}(-1)^{k+1}e^{-k} & \textrm{if $n$ is even,}
\end{array}
\right.
$$
and one then gets
$$
\limn\frac{\log(U(x_n))}{\log(x_n)}
=\left\{
\begin{array}{ll}
\displaystyle \frac{1}{1+e^{-1}} & \textrm{if $n$ is odd} \\
 & \\
\displaystyle \frac{e^{-1}}{1+e^{-1}} & \textrm{if $n$ is even,}
\end{array}
\right.
$$
which implies
$\nu(U)-\mu(U)\geq(1-e^{-1})\big/(1+e^{-1})>0$, hence the limit $\displaystyle \limx\frac{\log(U(x))}{\log(x)}$ does not exist and thus, applying Theorem CK, $U\not\in\M$.
\end{exm}

Now we give another relationship between \emph{O-RV} and $\M$.

\begin{prop}\label{teo:20141220:002}
Let $U:\Rset^+\to\Rset^+$ be a measurable function. 
If $U\in\textrm{O-RV}$ and the limit
$$
\displaystyle 
\limx\frac{\log(U(x))}{\log(x)}
$$
exists, then $U\in\M$.
\end{prop}

The relationships of $\M_\infty$ and $\M_{-\infty}$ with \emph{O-RV} are simpler.

\begin{prop}\label{teo:20150102:010}
For $\lambda\in\big\{\infty,-\infty\big\}$, \ $\M_\lambda\bigcap\textrm{O-RV}=\emptyset$.
\end{prop}


\section{Proofs}
\label{SectionProofs}

\begin{proof}[\textbf{Proof of Theorem \ref{teo:20141211:001}}]~

\begin{itemize}
\item\emph{Proof of the necessary condition of (i)}~

Assume $U\in\M$ with $\rho_U=-\tau$.
Let $r\in\Rset$ such that $r\neq0$.

\begin{itemize}
\item\emph{If $r<\tau$}

Let $0<\epsilon<\tau-r$ and $\delta>0$.
By hypothesis, there exists a constant $x_0>1$ such that, for $x\geq x_0$, $\displaystyle U(x)\leq\delta x^{-\tau+\epsilon}$, and there exists
$x_1>1$ such that, for $x\geq x_1$, $\displaystyle U(x)\geq x^{-\tau-\epsilon}\big/\delta$.
Hence, setting $x_a:=\max(x_0,x_1)$, for $x\geq x_a$ and $t>1$,
$$
\displaystyle t^r\frac{U(tx)}{U(x)}\leq \delta^2t^r(tx)^{-\tau+\epsilon}x^{\tau+\epsilon}=\delta^2t^{-\tau+r+\epsilon}x^{2\epsilon}\textrm{,}
$$
and the assertion then follows as $t\to\infty$ since $-\tau+r+\epsilon<0$.

\item\emph{If $r>\tau$}

Let $0<\epsilon<r-\tau$ and $\delta>0$.
By hypothesis, there exists a constant $x_0>1$ such that, for $x\geq x_0$, $\displaystyle U(x)\leq\delta x^{-\tau+\epsilon}$, and there exists
$x_1>1$ such that, for $x\geq x_1$, $\displaystyle U(x)\geq x^{-\tau-\epsilon}\big/\delta$.
Hence, setting $x_a:=\max(x_0,x_1)$, for $x\geq x_a$ and $t>1$,
$$
\displaystyle t^r\frac{U(tx)}{U(x)}\geq \frac{1}{\delta^2}t^r(tx)^{-\tau-\epsilon}x^{\tau-\epsilon}=\frac{1}{\delta^2}t^{r-\tau-\epsilon}x^{-2\epsilon}\textrm{,}
$$
and the assertion then follows as $t\to\infty$ since $r-\tau-\epsilon>0$.

\end{itemize}

\item\emph{Proof of the sufficient condition of (i)}

Let $\delta>0$ and $\eta>0$.

One the one hand, since $\tau-\delta\big/2<\tau$, by hypothesis, there exists a constant $x_a>1$ such that, for $x\geq x_a$,
$
\displaystyle
\limt t^{\tau-\delta/2}\frac{U(xt)}{U(x)}=0
$.
Hence, given $x\geq x_a$, there exists $t_a=t_a(x)>1$ such that, for $t\geq t_a$, $t^{\tau-\delta/2}U(tx)\leq \eta U(x)$, or
\begin{equation}\label{eq:20141211:010}
\frac{U(tx)}{(tx)^{-\tau+\delta}}\leq\eta\frac{x^{\tau-\delta}U(x)}{t^{\delta/2}}\textrm{.}
\end{equation}

One the other hand, since $\tau+\delta\big/2>\tau$, by hypothesis, there exists a constant $x_b>1$ such that, for $x\geq x_b$,
$
\displaystyle
\limt t^{\tau+\delta/2}\frac{U(xt)}{U(x)}=\infty
$.
Hence, given $x\geq \max(x_a,x_b)$, there exists $t_b=t_b(x)>1$ such that, for $t\geq t_b$, $t^{\tau+\delta/2}U(tx)\geq \eta U(x)$, or
\begin{equation}\label{eq:20141211:011}
\frac{U(tx)}{(tx)^{-\tau-\delta}}\geq\eta x^{\tau+\delta}U(x)t^{\delta/2}\textrm{.}
\end{equation}
Combining \eqref{eq:20141211:010} and \eqref{eq:20141211:011}, given $x\geq \max(x_a,x_b)$ and for $t\geq\max(t_a,t_b)$, and using the change of variable $y=tx$ with $y\to\infty$ as $t\to\infty$, provide, for $\delta>0$,
$$
\limy\frac{U(y)}{y^{-\tau+\delta}}=0
\quad\textrm{and}\quad
\limy\frac{U(y)}{y^{-\tau-\delta}}=\infty\textrm{,}
$$
which implies that $U\in\M$ with $\rho_U=-\tau$.

\item\emph{Proof of the necessary condition of (ii)}

Let $r\in\Rset$ and $\eta>0$.
Set $r'<-r$.
Since $U\in\M_\infty$ there exists a constant $x_0>1$ such that, for $x\geq x_0$, $\displaystyle U(x)\leq\eta x^{r'}$.
Hence, for $t>1$,
$$
\displaystyle t^r\frac{U(tx)}{U(x)}\leq \eta\frac{t^{r+r'}x^{r'}}{U(x)}\textrm{,}
$$
and the assertion then follows as $t\to\infty$ since $r+r'<0$.

\item\emph{Proof of the sufficient condition of (ii)}

Let $r\in\Rset$.
Taking $r'<-r$, by hypothesis, there exists a constant $x_0>1$ such that, for $x\geq x_0$,
$
\displaystyle
\limt t^{r'}\frac{U(xt)}{U(x)}=0
$.
Hence, for $\eta>0$, there exists a constant $t_0>1$ such that, for $t\geq t_0$, $t^{r'}U(tx)\leq \eta U(x)$, or
$$
\frac{U(tx)}{(tx)^r}\leq\eta\frac{U(x)}{x^rt^{r+r'}}\textrm{.}
$$
Using the change of variable $y=tx$ and noting that $y\to\infty$ as $t\to\infty$ give, for $r\in\Rset$, being $r+r'>0$,
$$
\limy\frac{U(y)}{y^r}=0\textrm{,}
$$
which means that $U\in\M_\infty$.

\item\emph{Proof of the necessary condition of (iii)}

Let $r\in\Rset$ and $\eta>0$.
Set $r'>-r$.
Since $U\in\M_{-\infty}$ there exists a constant $x_0>1$ such that, for $x\geq x_0$, $\displaystyle U(x)\geq\eta x^{r'}$.
Hence, for $t>1$,
$$
\displaystyle t^r\frac{U(tx)}{U(x)}\geq \eta\frac{x^{r'}}{U(x)}t^{r+r'}\textrm{,}
$$
and the assertion then follows as $t\to\infty$ since $r+r'>0$.

\item\emph{Proof of the sufficient condition of (iii)}

Let $r\in\Rset$.
Taking $r'<-r$, by hypothesis, there exists a constant $x_0>1$ such that, for $x\geq x_0$,
$
\displaystyle
\limt t^{r'}\frac{U(xt)}{U(x)}=\infty
$.
Hence, for $\eta>0$, there exists a constant $t_0>1$ such that, for $t\geq t_0$, $t^{r'}U(tx)\geq \eta U(x)$, or
$$
\frac{U(tx)}{(tx)^r}\geq\eta\frac{U(x)}{x^r}t^{-r-r'}\textrm{.}
$$
Using the change of variable $y=tx$ and noting that $y\to\infty$ as $t\to\infty$ give, for $r\in\Rset$, being $-r-r'>0$,
$$
\limy\frac{U(y)}{y^r}=0\textrm{,}
$$
which means that $U\in\M_\infty$.
\end{itemize}
\end{proof}

\begin{proof}[\textbf{Proof of Corollary \ref{cor:20141225:001}}]~

Let $U\in\textrm{RV}$ with tail index $\rho$.
Then, for $t>1$,
$$
\limx t^r\frac{U(tx)}{U(x)}=t^{r+\rho}\textrm{,}
$$
which implies that, for $\epsilon>0$, there exists a constant $x_0>1$ such that, for $x\geq x_0$,
$$
t^{r+\rho}-\epsilon\leq t^r\frac{U(tx)}{U(x)}\leq t^{r+\rho}+\epsilon\textrm{.}
$$
Hence, setting $\tau=-\rho$, gives, on the one hand, for $r<\tau$,
$$
-\epsilon\leq\limt t^r\frac{U(tx)}{U(x)}\leq\epsilon\textrm{,}
$$
which implies $\displaystyle \limt t^r\frac{U(tx)}{U(x)}=0$ taking $\epsilon$ arbitrary, and, on the other hand, for $r>\tau$,
$$
\limt t^r\frac{U(tx)}{U(x)}=\infty\textrm{.}
$$
Therefore one has, applying Theorem \ref{teo:20141211:001}, that $U\in\M$ with $\rho_U=\rho$.

Finally, a function belonging to $\M$ but not to \emph{RV} is for instance the function given in Example \ref{exm:20150103:001}.
\end{proof}

\begin{proof}[\textbf{Proof of Proposition \ref{prop:20141211:001}}]~

\begin{itemize}
\item\emph{Proof of (i)}

Let $\epsilon>0$. By definition of $U\in\M$ with $\rho_U=-\tau$, there exist constants $x_a,x_b>1$ such that,
$$
\textrm{for $x\geq x_a$, } U(x)\leq x^{-\tau+\epsilon}\textrm{,}
\quad\textrm{and,}\quad
\textrm{for $x\geq x_b$, } U(x)\geq x^{-\tau-\epsilon}\textrm{.}
$$
So, for $x\geq x_0:=\max(x_a,x_b)$, $x^{-\tau-\epsilon}\leq U(x)\leq x^{-\tau+\epsilon}$.
Hence, for any $x_0\leq c< d<\infty$, one has, setting $M_c:=\min(c^{-\tau-\epsilon},d^{-\tau+\epsilon})$ and $M_d:=\max(c^{-\tau-\epsilon},d^{-\tau+\epsilon})$, that $U$ satisfies $M_c\leq U(x)\leq M_d$ for any $x\in[c;d]$.

\item\emph{Proof of (ii)}

Let $\epsilon>0$. By definition of $U\in\M_\infty$, there exists a constant $x_0>1$ such that, for $x\geq x_0$,
$
U(x)\leq x^{\epsilon}\textrm{.}
$
Hence, for any $c\geq x_0$, one has, setting $M_c:=c^{\epsilon}$, that $U$ satisfies $U(x)\leq M_c$ for any $x\in[c;\infty)$.

\item\emph{Proof of (iii)}

Let $\epsilon>0$. By definition of $U\in\M_{-\infty}$, there exists a constant $x_0>1$ such that, for $x\geq x_0$,
$
U(x)\geq x^{\epsilon}\textrm{.}
$
Hence, for any $d\geq x_0$, one has, setting $M_d:=d^{\epsilon}$, that $U$ satisfies $U(x)\geq M_d$ for any $x\in[d;\infty)$.
\end{itemize}
\end{proof}

\begin{proof}[\textbf{Proof of Theorem \ref{teo:20141211:002}}]~

Let $\mu$ be the Lebesgue measure on $\Rset$.

\begin{itemize}
\item\emph{Proof of (i)}

Let $U\in\M$ with $\rho_U=-\tau$ and let $r<\tau$.
Applying Theorem \ref{teo:20141211:001}, (i), there exists $x_a>1$ such that, for $x\geq x_a$,
$$
\limt t^r\frac{U(tx)}{U(x)}=0\textrm{.}
$$
Let $x_a\leq c< d<\infty$.
Then using Egoroff's theorem (see e.g. \cite{Billingsley}), there exists a measurable $A\subseteq[c;d]$ of a positive measure such that
$$
\limt\sup_{x\in A}t^r\frac{U(tx)}{U(x)}=0\textrm{.}
$$
Let us prove by contradiction that the previous limit holds on $[c;d]$.
Then suppose that there exist $\epsilon>0$, $\big\{x_n\big\}_{n\in\Nset}\subseteq[c;d]$, and $\big\{t_n\big\}_{n\in\Nset}\subseteq\Rset^+$ such that $t_n\to\infty$ and
\begin{equation}\label{eq:20141226:31}
\limn t_n^r\frac{U(t_nx_n)}{U(x_n)}>\epsilon\textrm{.}
\end{equation}
By Proposition \ref{propArandelovic2004} one has, denoting $\log(A)=\big\{\log(x):x\in A\big\}$ and noting that $\log(A)$ has a positive measure,
$$
\mu\left(\limsupn(\log(A)-\log(x_n))\right)\geq\mu\left(\log A\right)>0\textrm{,}
$$
which implies that there exist a constant $\log(u)\in\Rset$ and a subsequence $\big\{x_{n_i}\big\}_{i\in\Nset}\subseteq\big\{x_n\big\}_{n\in\Nset}$ such that $\log(x_{n_i})+\log(u)\in \log(A)$,
i.e. $u\,x_{n_i}\in A$.
Note that $u>0$.

By Proposition \ref{prop:20141211:001}, (i), 
there exist $0<M_c\leq M_d<\infty$ such that $M_c\leq U(x)\leq M_d$, $x\in(c;d)$.
Hence, one then has
$$
t_{n_i}^r\frac{U(t_{n_i}x_{n_i})}{U(x_{n_i})}=\left(\frac{t_{n_i}}{u}\right)^r\frac{U\left(\frac{t_{n_i}}{u}\,ux_{n_i}\right)}{U(ux_{n_i})}u^r\frac{U(ux_{n_i})}{U(x_{n_i})}
\leq\left(\frac{t_{n_i}}{u}\right)^r\frac{U\left(\frac{t_{n_i}}{u}\,ux_{n_i}\right)}{U(ux_{n_i})}u^r\frac{M_d}{M_c}\textrm{.}
$$
Noting that $\displaystyle \left(\frac{t_{n_i}}{u}\right)^r\frac{U((t_{n_i}\big/ u)\,ux_{n_i})}{U(ux_{n_i})}\to0$ since $u\,x_{n_i}\in A$ and $t_{n_i}\big/ u\to\infty$ as $n_i\to\infty$
provide
$
\displaystyle t_{n_i}^r\frac{U(t_{n_i}x_{n_i})}{U(x_{n_i})}\to0
$ as $n_i\to\infty$, which contradicts \eqref{eq:20141226:31}.

\item\emph{Proof of (ii)}

Let $U\in\M$ with $\rho_U=-\tau$ and let $r<\tau$.
Applying Theorem \ref{teo:20141211:001}, (i), there exists $x_b>1$ such that, for $x\geq x_b$,
$$
\limt t^r\frac{U(tx)}{U(x)}=\infty\textrm{.}
$$
Let $x_b\leq c< d<\infty$ and let $\big\{\epsilon_m\big\}_{m\in\Nset}$ be a strictly increasing sequence of positive numbers such that $\epsilon_m\to\infty$ as $m\to\infty$.
Then using Egoroff's theorem, there exists a measurable $A_m\subseteq[c;d]$, $m\in\Nset$, of a positive measure such that
$$
\limt\inf_{x\in A_m}t^r\frac{U(tx)}{U(x)}\geq\epsilon_m\textrm{.}
$$
Let us prove
$$
\limt\inf_{x\in[c;d]}t^r\frac{U(tx)}{U(x)}=\infty\textrm{}
$$
by contradiction.
Then suppose that there exist $\delta>0$, $\big\{x_n\big\}_{n\in\Nset}\subseteq[c;d]$, and $\big\{t_n\big\}_{n\in\Nset}\subseteq\Rset^+$ such that $t_n\to\infty$ and
\begin{equation}\label{eq:20141226:32}
\limn t_n^r\frac{U(t_nx_n)}{U(x_n)}<\delta\textrm{.}
\end{equation}
By Proposition \ref{propArandelovic2004} one has, denoting $\log(A_m)=\big\{\log(x):x\in A_m\big\}$, $m\in\Nset$, and noting that $\log(A_m)$ has a positive measure,
$$
\mu\left(\limsupn(\log(A_m)-\log(x_n))\right)\geq\mu\left(\log A_m\right)>0\textrm{,}
$$
which implies, for $m\in\Nset$, that there exist a constant $\log(u_m)\in\Rset$ and a subsequence $\big\{x_{n_{m,i}}\big\}_{i\in\Nset}\subseteq\big\{x_n\big\}_{n\in\Nset}$ such that $\log(x_{n_{m,i}})+\log(u_m)\in \log(A_m)$,
i.e. $u_m\,x_{n_{m,i}}\in A_m$.
Note that $u_m>0$ and $c\big/d\leq u_m\leq d\big/c$, $m\in\Nset$.

By Proposition \ref{prop:20141211:001}, (i), 
there exist $0<M_c\leq M_d<\infty$ such that $M_c\leq U(x)\leq M_d$ for $x\in(c;d)$.
Hence, one then has
$$
t_{n_{m,i}}^r\frac{U(t_{n_{m,i}}x_{n_{m,i}})}{U(x_{n_{m,i}})}=\left(\frac{t_{n_{m,i}}}{u_m}\right)^r\frac{U\left(\frac{t_{n_{m,i}}}{u_m}\,u_mx_{n_{m,i}}\right)}{U(u_mx_{n_{m,i}})}u_m^r\frac{U(u_mx_{n_{m,i}})}{U(x_{n_{m,i}})}
\geq\epsilon_m\left(\frac{c}{d}\right)^r\frac{M_c}{M_d}\textrm{,}
$$
implying 
$
\displaystyle t_{n_{m,i}}^r\frac{U(t_{n_{m,i}}x_{n_{m,i}})}{U(x_{n_{m,i}})}\to\infty
$ as $m\to\infty$, which contradicts \eqref{eq:20141226:32}.

\item\emph{Proof of (iii)}

Let $U\in\M_\infty$ and let $r\in\Rset$.
Applying Theorem \ref{teo:20141211:001}, (ii), there exists $x_0>1$ such that, for $x\geq x_0$,
$$
\limt t^r\frac{U(tx)}{U(x)}=0\textrm{.}
$$
Let $x_0\leq c< d<\infty$.

On the one hand, by hypothesis, there exists a constant $M_d>0$ such that, for $x\in[1;d]$, $U(x)\geq M_d$.
On the other hand, by Proposition \ref{prop:20141211:001}, (ii), there exists a constant $M_c>0$ such that, for $x\in[c;\infty)$, $U(x)\leq M_c$.
Combining these inequalities gives, for $x\in[c;d]$, $M_d\leq U(x)\leq M_c$.
Hence a proof similar to the one given to prove (i) can be done to conclude that $\displaystyle \limx t^r\sup_{x\in[c;d]}\frac{U(tx)}{U(x)}=0$.

\item\emph{Proof of (iv)}

Let $U\in\M_{-\infty}$ and let $r\in\Rset$.
Applying Theorem \ref{teo:20141211:001}, (iii), there exists $x_0>1$ such that, for $x\geq x_0$,
$$
\limt t^r\frac{U(tx)}{U(x)}=\infty\textrm{.}
$$
Let $x_0\leq c< d<\infty$.

On the one hand, by hypothesis, there exists a constant $M_d>0$ such that, for $x\in[1;d]$, $U(x)\leq M_d$.
On the other hand, by Proposition \ref{prop:20141211:001}, (iii), there exists a constant $M_c>0$ such that, for $x\in[c;\infty)$, $U(x)\geq M_c$.
Combining these inequalities gives, for $x\in[c;d]$, $M_c\leq U(x)\leq M_d$.
Hence a proof similar to the one given to prove (ii) can be done to conclude that $\displaystyle \limx t^r\inf_{x\in[c;d]}\frac{U(tx)}{U(x)}=\infty$.

\end{itemize}
\end{proof}

\begin{proof}[\textbf{Proof of Theorem CK}]~

Let $U:\Rset^+\to\Rset^+$ be a measurable function.
\begin{itemize}
\item\emph{Proof of (i) $\Rightarrow$ (ii)}

Let $\epsilon>0$ and $U\in\M$ with $\rho_U=\tau$.
One has, by definition, that
$$
\limx\frac{U(x)}{x^{\rho+\epsilon}}=0
\quad\textrm{and}\quad
\limx\frac{U(x)}{x^{\rho-\epsilon}}=\infty\textrm{.}
$$
Hence, there exists $x_0\geq1$ such that, for $x\geq x_0$,
$$
U(x)\leq\epsilon x^{\tau+\epsilon}
\quad\textrm{and}\quad
U(x)\geq\frac{1}{\epsilon}x^{\tau-\epsilon}\textrm{.}
$$
Applying the logarithm function to these inequalities and dividing them by $\log(x)$ (with $x>1$) provide
$$
\frac{\log\left(U(x)\right)}{\log(x)}\leq\frac{\log\left(\epsilon\right)}{\log(x)}+\tau+\epsilon
\quad\textrm{and}\quad
\frac{\log\left(U(x)\right)}{\log(x)}\geq-\frac{\log\left(\epsilon\right)}{\log(x)}+\tau-\epsilon\textrm{,}
$$
and, one then has
$$
\limsupx\frac{\log\left(U(x)\right)}{\log(x)}\leq\tau+\epsilon
\quad\textrm{and}\quad
\liminfx\frac{\log\left(U(x)\right)}{\log(x)}\geq\tau-\epsilon\textrm{,}
$$
from which one gets, taking $\epsilon$ arbitrary,
$$
\tau\leq\liminfx\frac{\log\left(U(x)\right)}{\log(x)}\leq\limsupx\frac{\log\left(U(x)\right)}{\log(x)}\leq\tau\textrm{,}
$$
and the assertion follows.

\item\emph{Proof of (ii) $\Rightarrow$ (iii)}

Let $0<\epsilon<1\big/2$.
Assume $U$ satisfies
$
\displaystyle \limx\frac{\log\left(U(x)\right)}{\log(x)}=\tau\textrm{.}
$
Let $\gamma$ a measurable function with support $\Rset^+$ such that $\gamma(x)\to0$ as $x\to\infty$, and let $b>1$.
Applying \thelr\ to the ratio gives
$$
\limx\left(\gamma(x)+\frac{\int_b^x\frac{\log(U(s))}{\log(s)}\frac{ds}{s}}{\log(x)}\right)
=\limx\frac{\log(U(x))}{\log(x)}=\tau
$$
First, suppose $\tau\neq0$, then
$$
\limx\frac{\log(U(x))}{\gamma(x)\log(x)+\int_b^x\frac{\log(U(s))}{\log(s)}\frac{ds}{s}}=1\textrm{,}
$$
and there exists $x_0>1$ such that, for $x\geq x_0$,
$$
\delta_U(x):=\frac{\log(U(x))}{\gamma(x)\log(x)+\int_b^x\frac{\log(U(s))}{\log(s)}\frac{ds}{s}}\geq1-\epsilon>0\textrm{.}
$$
Setting $x_1:=\max(b,x_0)$ and defining the functions, for $x\geq x_1$, $\alpha_U(x):=\gamma(x)\delta_U(x)\log(x)$ and $\beta_U(x):=\log(U(x))\big/\log(x)$, the assertion follows.

Now, suppose $\tau=0$.
Define the function $V(x):=xU(x)$, $x>0$, which clearly satisfies $\displaystyle \limx\frac{\log(V(x))}{\log(x)}=1\neq0$.
Hence, applying to $V$ the previous proof for $U$ when $\tau\neq0$ gives that there exist $x_{1,V}\geq b_V>1$ and measurable functions $\alpha_V$, $\beta_V$, and $\delta_V$ satisfying, as $x\to\infty$,
$$
\alpha_V(x)\big/\log(x)\to0\textrm{,}\quad\beta_V(x)\to1\textrm{,}\quad\delta_V(x)\to1\textrm{,}
$$
such that, for $x\geq x_{1,V}$,
$$
V(x)=\exp\left\{\alpha_V(x)+\delta_V(x)\int_{b_V}^x\beta_V(s)\frac{ds}{s}\right\}\textrm{.}
$$
Defining, when $\tau=0$, the constant $x_{1,U}:=x_{1,V}$ and the functions $\alpha_U(x):=\alpha_V(x)+\log(x)\big(\delta_V(x)-1\big)$, $\beta_U(x):=\beta_V(x)-1$, and $\delta_U(x):=\delta_V(x)$, the assertion follows.

\item\emph{Proof of (iii) $\Rightarrow$ (i)}

Suppose there exist $b>1$ and measurable functions $\alpha$, $\beta$, and $\delta$ satisfying, as $x\to\infty$,
$$
\alpha(x)\big/\log(x)\to0\textrm{,}\quad\beta(x)\to\tau\textrm{,}\quad\delta(x)\to1\textrm{,}
$$
such that
$$
U(x)=\exp\left\{\alpha(x)+\delta(x)\int_b^x\beta(s)\frac{ds}{s}\right\}\textrm{,\quad$x\geq x_1$ for some $x_1\geq b$.}
$$
Let $\epsilon>0$ sufficiently small such that $2\epsilon\big(\tau+\epsilon\big/4\big)\leq1$ and $2\epsilon\big(\tau-\epsilon\big/4\big)\geq-1$.
Then there exist $x_a>1$ such that, for $x\geq x_a$, $\big|\alpha(x)\big/\log(x)\big|\leq\epsilon\big/4$, $x_b>1$ such that, for $x\geq x_b$, $\big|\beta(x)-\tau\big|\leq\epsilon\big/4$, 
and $x_c>1$ such that, for $x\geq x_c$, $\big|\delta(x)-1\big|\leq\epsilon^2/4$.

On the one hand, writing, for $x\geq x_0:=\max(b,x_a,x_b,x_c)$,
\begin{eqnarray*}
\lefteqn{\frac{U(x)}{x^{\tau+\epsilon}}=\exp\left\{-(\tau+\epsilon)\log(x)+\alpha(x)+\delta(x)\int_b^x\beta(s)\frac{ds}{s}\right\}} \\
 & & =\exp\left\{\log(x)\left(\frac{\alpha(x)}{\log(x)}-\frac{\epsilon}{2}\right)+\delta(x)\int_b^x\beta(s)\frac{ds}{s}-\left(\tau+\frac{\epsilon}{2}\right)\log(x)\right\} \\
 & & \leq\exp\left\{-\frac{\epsilon}{4}\log(x)+\delta(x)\int_b^{x_0}\beta(s)\frac{ds}{s}+\delta(x)\left(\tau+\frac{\epsilon}{4}\right)\left(\log(x)-\log(x_0)\right)-\left(\tau+\frac{\epsilon}{2}\right)\log(x)\right\}
\end{eqnarray*}
and noting that
$$
\delta(x)\left(\tau+\frac{\epsilon}{4}\right)-\left(\tau+\frac{\epsilon}{2}\right)=
\big(\delta(x)-1\big)\left(\tau+\frac{\epsilon}{4}\right)-\frac{\epsilon}{4}\leq-\frac{\epsilon}{8}
$$
give
\begin{equation}\label{eq:20150106:01}
\limx\frac{U(x)}{x^{\tau+\epsilon}}=0\textrm{.}
\end{equation}
On the one hand, writing, for $x\geq x_0:=\max(b,x_a,x_b,x_c)$,
\begin{eqnarray*}
\lefteqn{\frac{U(x)}{x^{\tau-\epsilon}}=\exp\left\{-(\tau-\epsilon)\log(x)+\alpha(x)+\delta(x)\int_b^x\beta(s)\frac{ds}{s}\right\}} \\
 & & =\exp\left\{\log(x)\left(\frac{\alpha(x)}{\log(x)}+\frac{\epsilon}{2}\right)+\delta(x)\int_b^x\beta(s)\frac{ds}{s}-\left(\tau-\frac{\epsilon}{2}\right)\log(x)\right\} \\
 & & \geq\exp\left\{\frac{\epsilon}{4}\log(x)+\delta(x)\int_b^{x_0}\beta(s)\frac{ds}{s}+\delta(x)\left(\tau-\frac{\epsilon}{4}\right)\left(\log(x)-\log(x_0)\right)-\left(\tau-\frac{\epsilon}{2}\right)\log(x)\right\}
\end{eqnarray*}
and noting that
$$
\delta(x)\left(\tau-\frac{\epsilon}{4}\right)-\left(\tau-\frac{\epsilon}{2}\right)=
\big(\delta(x)-1\big)\left(\tau-\frac{\epsilon}{4}\right)+\frac{\epsilon}{4}\geq\frac{\epsilon}{8}
$$
give
\begin{equation}\label{eq:20150106:02}
\limx\frac{U(x)}{x^{\tau-\epsilon}}=\infty\textrm{.}
\end{equation}
Combining \eqref{eq:20150106:01} and \eqref{eq:20150106:02} provides $U\in\M$ with $\rho_U=\tau$.
\end{itemize}
\end{proof}

\begin{proof}[\textbf{Proof of Proposition \ref{teo:20141220:002}}]~

Let $U:\Rset^+\to\Rset^+$ be a measurable function.

Assume $U\in\textrm{\emph{O-RV}}$ and the limit
$
\displaystyle 
\limx\frac{\log(U(x))}{\log(x)}
$
exists.
Applying Theorem CK gives $U\in\M$ with $\displaystyle \rho_U=\limx\frac{\log(U(x))}{\log(x)}$.
\end{proof}

\begin{proof}[\textbf{Proof of Proposition \ref{teo:20150102:010}}]
We will prove the proposition for $\lambda=\infty$. The proof for $\lambda=-\infty$ is similar.

Let us prove it by contradiction.
Assume there exists $U\in\M_\infty\bigcap\textrm{\emph{O-RV}}$.

By assumption $U\in\M$, we have, for $\rho\in\Rset$ and $\delta>0$, there exists $x_0>1$ such that, for $x\geq x_0$, $U(x)\leq cx^\rho$.
Applying the logarithm function to this inequality, dividing it by $\log(x)$, $x>1$, and taking the limit $x\to\infty$ give
$$
\limx\frac{\log(U(x))}{\log(x)}\leq\rho\textrm{.}
$$
Taking $\rho$ arbitrary provides
\begin{equation}\label{eq:20150105:01}
\limx\frac{\log(U(x))}{\log(x)}=-\infty\textrm{.}
\end{equation}
Now, by assumption $U\in\textrm{\emph{O-RV}}$, applying Proposition \ref{prop:20141220:000}, (i) $\Rightarrow$ (ii), there exist $\alpha,\beta\in\Rset$ and $x_1>1$, $c>0$ such that, for all $t\geq1$ and $x\geq x_1$,
$$
c^{-1}t^\beta\leq\frac{U(tx)}{U(x)}\leq ct^\alpha\textrm{.}
$$
Hence applying to these inequalities the logarithm function, dividing them by $\log(t)$, $t>0$, and taking the limit $t\to\infty$ give
$$
\left|\limt\frac{\log(U(t))}{\log(t)}\right|\leq\max\big\{|\alpha|,|\beta|\big\}<\infty\textrm{,}
$$
which contradicts \eqref{eq:20150105:01}.
The proposition is proved.
\end{proof}

\section{Conclusion}

A new characterization of the class $\M$ introduced in \cite{NN2014}, a strict larger class than the class of regularly varying functions (RV), was proved, and it was extended to the classes $\M_\infty$ and $\M_{-\infty}$.
This characterization together with other two given by Cadena and Kratz in \cite{NN2014} allowed the study of relationships between $\M$ and the well-known class \emph{O-RV}, another extension of RV.
It was found that these classes satisfy $\M$ \ $\not\subseteq$ \ \emph{O-RV} and \emph{O-RV} \ $\not\subseteq$ \ $\M$, and
necessary conditions to have inclusions were provided.
Relationships among \emph{O-RV} and $\M_\infty$ and $\M_{-\infty}$ were provided.

Note that any result obtained here can be applied to positive and measurable functions with finite support 
by using the change of variable $y =1\big/(x^*_U-x)$ for $x < x^*_U$ where $x^*_U$ is the \emph{endpoint} of $U$ defined by $x^*_U:=\sup\big\{x:U(x)>0\big\}$.

\section*{Acknowledgments} 
The author gratefully acknowledges the support of SWISS LIFE through its ESSEC research program on 'Consequences of the population ageing on the insurances loss'.




\end{document}